\documentclass[12pt]{article}
\usepackage{amssymb, amsmath, latexsym, a4}
\usepackage[latin1]{inputenc}
\usepackage[all]{xy}
 \usepackage[T1]{fontenc}
\usepackage[usenames]{color}
\begin{document}

\newcommand{\rdg}{\hfill $\Box $}

\newtheorem{definition}{Definition}[section]
\newtheorem{theorem}[definition]{Theorem}
\newtheorem{proposition}[definition]{Proposition}
\newtheorem{lemma}[definition]{Lemma}
\newtheorem{corollary}[definition]{Corollary}
\newtheorem{remark}[definition]{Remark}
\newtheorem{example}[definition]{Example}
\newtheorem{exercise}[definition]{Exercises}
\newcommand{\tp}{\otimes}
\newcommand{\N}{\mathbb{N}}
\newcommand{\Z}{\mathbb{Z}}
\newcommand{\op}{\oplus}
\newcommand{\n}{\underline n}
\newcommand{\es}{{\frak S}}
\newcommand{\ef}{\frak F}
\newcommand{\qu}{\frak Q} \newcommand{\ga}{\frak g}
\newcommand{\la}{\lambda}
\newcommand{\ig}{\frak Y}
\newcommand{\te}{\frak T}
\newcommand{\cok}{{\sf Coker}}
\newcommand{\Hom}{{\sf Hom}}
\newcommand{\im}{{\sf Im}}
\newcommand{\ext}{{\sf Ext}}
\newcommand{\ho}{{\sf H_{AWB}}}
\newcommand{\HH}{{\sf Hoch}}
\newcommand{\adu}{{\rm AWB}^!}

\def\Im{\operatorname{Im}}
\def\Ker{\operatorname{Ker}}
\def\id{\operatorname{id}}
\def\Coker{\operatorname{Coker}}
\def\Der{\operatorname{Der}}
\def\hom{\operatorname{Hom}}
\newcommand{\K}{\mathbb{K}}
\newcommand{\ad}{\operatorname{ ad}}

\newcommand{\ele}{\cal L} \newcommand{\as}{\cal A} \newcommand{\ka}{\cal
K}\newcommand{\eme}{\cal M} \newcommand{\pe}{\cal P}

\newcommand{\pn}{\par \noindent}
\newcommand{\pbn}{\par \bigskip \noindent}
\bigskip\bigskip

\centerline {\Large {\bf  A non-abelian tensor product of Hom-Lie algebras}}

\

\centerline {\bf J. M. Casas$^{(1)}$, E. Khmaladze$^{(2)}$ and N. Pacheco Rego$^{(3)}$}

\bigskip \bigskip
\centerline{$^{(1)}$Dpto.  Matem\'atica Aplicada I, Univ. de Vigo, 36005 Pontevedra, Spain}
\centerline{e-mail address: \tt jmcasas@uvigo.es}
\medskip

\centerline{$^{(2)}$A. Razmadze Math. Inst. of I. Javakhishvili Tbilisi State University,}
\centerline{Tamarashvili Str. 6, 0177 Tbilisi, Georgia and}
\centerline{Dpto.  Matem\'atica Aplicada I, Univ. de Vigo, 36005 Pontevedra, Spain}
\centerline{e-mail address: \tt e.khmal@gmail.com}
\medskip

\centerline{$^{(3)}$IPCA, Dpto. de Ciências, Campus do IPCA,
 Lugar do Aldão}
\centerline{4750-810 Vila Frescainha, S. Martinho, Barcelos,
 Portugal}
\centerline{e-mail address: \tt natarego@gmail.com}

\bigskip \bigskip

\par

{\bf Abstract:} Non-abelian tensor product of Hom-Lie algebras is constructed and studied.  This tensor product is used to describe universal ($\alpha$-)central extensions of Hom-Lie algebras and to establish a relation between cyclic and Milnor cyclic homologies of Hom-associative algebras satisfying certain additional condition.

\bigskip


\bigskip
{\it A. M. S. Subject Class. (2010):} 17A30, 17B55, 17B60, 18G35, 18G60

\section*{Introduction}

The concept of a Hom-Lie algebra was initially introduced in \cite{HLS}  motivated by discretization of vector fields via twisted derivations.
A Hom-Lie algebra is a non-associative algebra satisfying the skew symmetry and the Jacobi identity twisted by a map. When this map is the identity map, then
the definition of a Lie algebra is recovered. Thus, it is natural to seek for possible generalizations of known theories from Lie to Hom-Lie algebras.
In this context, recently there have been several works dealing with the study of Hom-Lie structures (see \cite{CaInPa}, \cite{JL}, \cite{MS} - \cite{Yau2}).

In this paper we introduce a non-abelian tensor product of Hom-Lie algebras generalizing the non-abelian tensor product of Lie algebras \cite{El1} and investigate its properties. In particular, we study its relation to the low dimensional homology of Hom-Lie algebras developed in \cite{Yau1, Yau2}. We use this tensor product in the description of the universal ($\alpha$-) central extensions of Hom-Lie algebras considered in \cite{CaInPa}. We give an application of our non-abelian tensor product of Hom-Lie algebras to cyclic homology of Hom-associative algebras \cite{Yau1}. Namely, for Hom-associative algebras satisfying an additional condition, which we chose to call $\alpha$-identity condition, we establish a relation between cyclic and Milnor cyclic homologies in terms of exact sequences.

Note that not all classical results can be generalized from Lie to Hom-Lie algebras, for example, results on universal central extensions of Lie algebras cannot be extended directly to Hom-Lie algebras and they are divided between universal central and universal $\alpha$-central extensions of Hom-Lie algebras (see Section \ref{section}). Further, in order to obtain Hom-algebra version of Guin's result relating cyclic and Milnor cyclic homology of associative algebras \cite{Gu}, we need to consider a subclass of Hom-associative algebras defined by the $\alpha$-identity condition (see Definition \ref{alfa condition}).

\subsection*{Notations} Throughout this paper we fix $\mathbb{K}$ as a ground field. Vector spaces are considered over $\mathbb{K}$ and linear maps are $\mathbb{K}$-linear maps. We write $\otimes$ (resp. $\wedge$) for the tensor product $\otimes_\mathbb{K}$ (resp. exterior product $\wedge_\mathbb{K}$ ) over $\mathbb{K}$. For any vector space (resp. Hom-Lie algebra) $L$, a subspace (resp. an ideal) $L'$ and $x\in L$ we write $\overline{x}$ to denote the coset $x+L'$.

\section{Hom-Lie algebras}
In this section we review some terminology and recall notions used in the paper. We mainly follow \cite{HLS, JL, MS, Yau1}, although with some modifications.

\subsection{Basic definitions}

\begin{definition}\label{def}
A Hom-Lie algebra $(L, \alpha_L)$ is a non-associative algebra $L$ together with a linear
map $\alpha_L : L \to L$ satisfying
\begin{align*}
&[x,y] = - [y,x], & \text{ (skew-symmetry)}\\
& [\alpha_L(x),[y,z]]+[\alpha_L(z),[x,y]]+[\alpha_L(y),[z,x]]=0  & \text{(Hom-Jacobi identity)}
\end{align*}
for all $x, y, z \in L$, where $[-,-]$ denotes the product in $L$.
\end{definition}

In this paper we only consider (the so called multiplicative) Hom-Lie algebras $(L, \alpha_L)$ such that $\alpha_L$ preserves the product, i.e. $\alpha_L[x,y]=[\alpha_L(x), \alpha_L(y)]$ for all $x, y \in L$.

\begin{example}\label{ejemplo 1} \
\begin{enumerate}
\item[a)] Taking $\alpha_L = \id_L$, Definition \ref{def} gives us the definition of a Lie algebra. Hence any Lie algebra $L$ can be considered as a Hom-Lie algebra $(L,\id_L)$.
\item[b)]   Let $V$ be a vector space and $\alpha_V:V\to V$ a linear map, then the pair $(V , \alpha_V)$ is called Hom-vector space. A Hom-vector space $(V,\alpha_V)$ together with the trivial product $[-,-]$ (i.e. $[x,y] = 0$ for any $x,y \in V$) is a Hom-Lie algebra $(V, \alpha_V)$ and it is called abelian Hom-Lie algebra.
\item[c)] Let $L$ be a Lie algebra, $[-,-]$ be the product in $L$  and $\alpha:L \to L$ be a Lie algebra endomorphism. Define $[-,-]_{\alpha} : L \otimes L \to L$ by $[x,y]_{\alpha} = \alpha_L[x,y]$, for all $x, y \in L$. Then $(L, \alpha)$ with the product $[-,-]_{\alpha}$ is a Hom-Lie algebra \cite[Theorem 5.3]{Yau1}.
\item[d)] Any Hom-associative algebra \cite{MS} becomes a  Hom-Lie algebra (see Section \ref{section 6} below).

\end{enumerate}
\end{example}

Hom-Lie algebras form a category ${\tt HomLie}$ whose {\it{morphisms}} from $(L,\alpha_L)$ to $(L',\alpha_{L'})$ are algebra homomorphisms $f : L \to L'$ such that $f \circ \alpha_L = \alpha_{L'} \circ f$. Clearly there is a full embedding ${\tt Lie} \hookrightarrow{\tt HomLie}$, $L\mapsto (L,\id_L)$, where  ${\tt Lie}$ denotes the category of Lie algebras.

It is a routine task to check that  ${\tt HomLie}$ satisfies the axioms of a semi-abelian category \cite{BB}. Consequently, the well-known Snake Lemma is valid for Hom-Lie algebras and we will use it in the sequel. Below we give the ad-hoc definitions of ideal, center, commutator, action and semi-direct product of Hom-Lie algebras and of course these notions agree with the respective general notions in the context of semi-abelian categories (see e.g. \cite{BoJaKe}).

\begin{definition}
A  Hom-Lie subalgebra $(H,\alpha_{H})$ of a Hom-Lie algebra $(L, \alpha_L)$ is a vector subspace $H$ of $L$ closed under the product, that is $[x,y] \in H$ for all $x, y \in H$, together with the endomorphism $\alpha_{H}:H\to H$ being the restriction of $\alpha_{L }$ on $H$. In such a case we may write $\alpha_{L\mid}$ for $\alpha_{H}$.

A  Hom-Lie subalgebra $(H,\alpha_{L\mid})$ of $(L, \alpha_L)$ is said to be an  ideal if $[x, y] \in H$ for any $x \in H$, $y\in L$.

If $(H,\alpha_{L\mid})$ is an ideal of a Hom-Lie algebra $(L,\alpha_L)$, then  $(L/H, \overline{\alpha}_L)$, where $\overline{\alpha }_L: L/H \to L/H$ is induced by $\alpha_L$, naturally inherits a structure of Hom-Lie algebra and it is called quotient Hom-Lie algebra.

Let $(H,\alpha_{L\mid})$  and $(K,\alpha_{L\mid})$  be  ideals of a Hom-Lie algebra $(L,\alpha_L)$. The commutator (resp. sum) of $(H,\alpha_{L\mid})$  and $(K,\alpha_{L\mid})$, denoted by $([H,K], \alpha_{{L\mid}})$ (resp. $(H+K, \alpha_{{L\mid}})$ ), is the  Hom-Lie subalgebra of $(L,\alpha_L)$ spanned by the elements $[h,k]$ (resp. $h+k$),  $h \in H$, $k \in K$.
\end{definition}

The following lemma is an easy exercise.

\begin{lemma}\label{ideales} Let $(H,\alpha_{L\mid})$ and $(K,\alpha_{L\mid})$ be  ideals of a Hom-Lie algebra
$(L,\alpha_L)$. The following statements hold:

\begin{enumerate} \label{ideales}

\item[a)] $(H \cap K, \alpha_{{L\mid}})$ and $(H+K, \alpha_{{L\mid}})$ are  ideals of $(L,\alpha_L)$;

\item[b)]  $[H,K] \subseteq H \cap K$;

\item[c)] If $\alpha_L$ is surjective, then $([H,K],\alpha_{{L\mid}})$ is an ideal of $(L,\alpha_L)$;

\item[d)] $([H,K],\alpha_{{L\mid}})$ is an ideal of $(H,\alpha_{L\mid})$  and $(K,\alpha_{L\mid})$. In particular, $([L,L],\alpha_{{L\mid}})$ is an ideal of $(L,\alpha_L)$;

\item[e)] $(\alpha_L(L),\alpha_{L{\mid}})$ is a Hom-Lie subalgebra of $(L,\alpha_L)$;

\item[f)] If $H, K \subseteq \alpha_L(L)$, then $([H,K],\alpha_{{L\mid}})$ is an ideal of $(\alpha_L(L),\alpha_{{L\mid}})$.
\end{enumerate}
\end{lemma}

\begin{definition}
The center of a  Hom-Lie algebra  $(L,\alpha_L)$ is the vector subspace
\[
Z(L) = \{ x \in L \mid [x, y] =0 \ \text{ for all}\ y \in L\}.
\]
\end{definition}

\begin{remark}
When $\alpha_L : L \to L$ is a surjective endomorphism, then $(Z(L),\alpha_{{L\mid}} )$ is an abelian Hom-Lie algebra and an ideal of $(L, \alpha_L)$.
\end{remark}

\subsection{Hom-action and semi-direct product}

 \begin{definition} \label{action}
 Let $(L,\alpha_L)$ and $(M, \alpha_M)$ be  Hom-Lie algebras. A Hom-action of $(L,\alpha_L)$ on  $(M, \alpha_M)$ is a linear map $L \otimes M \to M,$ $x \otimes m\mapsto {}^xm$, satisfying the following properties:
 \begin{enumerate}
 \item[a)] ${}^{[x,y]} \alpha_M(m) = {}^{\alpha_L(x)}({}^y m) - {}^{\alpha_L(y)} ({}^x m)$,
 \item [b)] ${}^{\alpha_L(x)} [m,m'] = [{}^x m, \alpha_M(m')]+[\alpha_M (m), {}^x m']$,
 \item [c)] $\alpha_M({}^x m) = {}^{\alpha_L(x)} \alpha_M(m)$
 \end{enumerate}
 for all $x, y \in L$ and $m, m' \in M$.

 The Hom-action is called trivial if ${}^xm=0$ for all $x\in L$ and $m\in M$.
 \end{definition}

 \begin{remark}
 If $(M, \alpha_M)$ is an abelian Hom-Lie algebra enriched with a Hom-action of $(L,\alpha_L)$, then  $(M, \alpha_M)$ is nothing else but a Hom-module over $(L,\alpha_L)$ (see \cite{Yau1} for the definition).
 \end{remark}

\begin{example}\label{ejemplo 2} \
\begin{enumerate}
\item[a)] Let $(L, \alpha_{L})$ be a Hom-subalgebra of a  Hom-Lie algebra $(K,\alpha_{K})$ and $(H,\alpha_{H})$ an ideal of $(K,\alpha_{K})$. Then there exists a Hom-action of $(L,\alpha_{L})$ on $(H,\alpha_{H})$ given by the product in $K$. In particular, there is a Hom-action of $(L,\alpha_{L})$ on itself given by the product in $L$.

\item[b)] Let $L$ and $M$ be Lie algebras. Any Lie action of $L$ on $M$ (see e.g. \cite{El1}) defines a Hom-action of $(L, \id_{L})$ on $(M, \id_{M})$.

    \item[c)] Let $L$ be a Lie algebra and $\alpha : L \to L$ be an endomorphism. Suppose $M$ is an $L$-module in the usual sense and the action of $L$ on $M$ satisfies the condition ${}^{\alpha(x)}m = {}^x m$, for all $x \in L$, $m \in M$. Then $(M, \id_M)$ is a Hom-module over the Hom-Lie algebra $(L, \alpha)$ considered in Example \ref{ejemplo 1} c).

    \noindent As an example of such $L$, $\alpha$ and $M$ we can consider $L$ to be the $2$-dimensional vector space  with basis $\{ e_1,e_2 \}$, together with the product $[e_1,e_2]=-[e_2,e_1]=e_1$ and zero elsewhere, $\alpha$ to be represented by the matrix $\left( \begin{array}{cc} 1 & 1 \\ 0 & 1 \end{array} \right)$, and $M$ to be the ideal of $L$ generated by $\{e_1\}$.

    \item [d)] Any homomorphism of Hom-Lie algebras $(L,\alpha_L)\to (M,\alpha_M)$ induces a Hom-action of $(L,\alpha_L)$ on $(M,\alpha_M)$ in
 the standard way by taking images of elements of $L$ and product in $M$.

 \item [e)] Let
  $0 \to (M,\alpha_M) \stackrel{i}\to (K,\alpha_K) \stackrel{\pi}\to (L,\alpha_L) \to 0$ be a split short exact sequence of Hom-Lie algebras, that is, there exists a homomorphism of Hom-Lie algebras $s:(L,\alpha_L)\to (K,\alpha_K)$ such that $\pi \circ s=\id_L$. Then there is a Hom-action of $(L,\alpha_L)$ on $(M,\alpha_M)$ defined in the standard way:  ${}^{x}m=i^{-1}[s(x),i(m)]$, $x\in L$, $m\in M$.
\end{enumerate}
\end{example}

\begin{definition}
 Given a Hom-action of a Hom-Lie algebra $(L,\alpha_L)$ on a Hom-Lie algebra $(M,\alpha_M)$ we define the semi-direct product Hom-Lie algebra, $(M \rtimes L, {\alpha_{\rtimes}})$, with the underlying  vector space $M \oplus L$, endowed with the product
\[
[(m_1,x_1),(m_2,x_2)]= ([m_1,m_2]+{^{\alpha_L(x_1)}}{m_2} - {^{\alpha_L(x_2)}}{m_1},[x_1,x_2]),
\]
together with the endomorphism ${\alpha_{\rtimes}} : M \rtimes L \to M \rtimes L$ given by ${\alpha_{\rtimes}} (m, x) = \left(\alpha_M(m), \alpha_L(x)\right)$ for all $x,x_1,x_2 \in L$ and $m, m_1, m_2 \in M$.
\end{definition}

Straightforward calculations show that $(M \rtimes L, {\alpha_{\rtimes}})$ indeed is a Hom-Lie algebra  and there is a short exact sequence of Hom-Lie algebras
\begin{equation}\label{semi-direct}
0 \to (M,\alpha_M) \stackrel{i}\to (M \rtimes L,\alpha_{\rtimes}) \stackrel{\pi}\to (L,\alpha_L) \to 0,
\end{equation}
where $i(m) = (m,0)$, $\pi(m,l)= l$. Moreover, $(M,\alpha_M)$ is an ideal of $(M \rtimes L, {\alpha_{\rtimes}})$ and this sequence splits by $s: (L,\alpha_L) \to (M \rtimes L,\alpha_{\rtimes})$, $s(l)=(0,l)$.
Then, as in  Example \ref{ejemplo 2} {\it e)}, the above sequence defines a Hom-action of $(  L,\alpha_{L})$  on $(M, \alpha_M)$ given by
\[
 {^l}m = i^{-1}\left[  \left(  0,l\right)  ,\left(  m,0\right)
\right]  =i^{-1}\left(  {^{\alpha_{L}\left(
l\right)}} m,0\right)  = {^{\alpha_{L}\left(  l\right)}} m.
\]
So, in general, the Hom-action of $\left(  L,\alpha_{L}\right)$ on $\left( M,\alpha_{M}\right)$ defined by the split short exact sequence (\ref{semi-direct}) does not coincide with the initial Hom-action of $(L,\alpha_{L})$ on $( M,\alpha_{M})$, but it coincides with the induced Hom-action of $(\alpha_{L}(L),\alpha_{L\mid})$ on  $( M,\alpha_{M})$.

\begin{definition}
Let $(M,\alpha_M)$ be a Hom-module over a Hom-Lie algebra $(L,\alpha_L)$. A derivation from $(L,\alpha_L)$ to $(M,\alpha_M)$ is a
linear map $d : L \to M$ satisfying
\begin{enumerate}
\item[a)] $d[x,y] = {^{\alpha_L(x)}}d(y) - {^{\alpha_L(y)}}d(x)$ for all $x, y \in L$,
\item [b)] $\alpha_M \circ d = d \circ \alpha_L$.
\end{enumerate}
We denote by $\Der_{\alpha}\left(L, M\right)$ the vector space of all
derivations from $(L, \alpha_{L})$ to $(M,\alpha_{M})$.
\end{definition}

The following lemma is straightforward.

 \begin{lemma}\label{derivation}
 Let $(M,\alpha_M)$ be a Hom-module over a Hom-Lie algebra $(L,\alpha_L)$. Then the projection $\theta : M \rtimes L \to M$, $\theta(m,l)=m$, is a derivation, where $(M,\alpha_M)$ is considered as a Hom-module over $(M \rtimes L,\alpha_{\rtimes})$ via the homomorphism  $\pi$ in (\ref{semi-direct}).
\end{lemma}

\begin{proposition} \label{correpresentabilidad}
Let $(M,\alpha_M)$ be a Hom-module over a Hom-Lie algebra $(L,\alpha_L)$. For every homomorphism of Hom-Lie algebras $f:(K,\alpha_K) \to (L,\alpha_L)$ and every derivation
$d:(K,\alpha_K) \to  (M,\alpha_M)$  there exists a unique homomorphism of  Hom-Lie algebras $h :  (K,\alpha_K) \to (M \rtimes
L,\alpha_{\rtimes})$ such that $\pi \circ h = f$ and $ \theta \circ h = d$, that is, the following diagram is commutative:
\[
\xymatrix{
& \left(  K,\alpha_{K} \right) \ar[dr]^{f}  \ar[d]^{h} \ar[ld]_{d} &          \\
 \left(  M,\alpha_{M}\right) \ar[r]^{i\ } & \ar@/^0.3pc/[l]^{\theta \ }     \left(  M \rtimes L,{\alpha_{\rtimes}}\right)  \ar[r]^{\ \pi} &  \left(  L,\alpha_{L}\right).    }\]
Here $(M,\alpha_M)$ is regarded as a Hom-module over $(K,\alpha_K)$ via $f$.

  Conversely, every homomorphism of  Hom-Lie algebras  $h:\left( K,\alpha_{K}\right)  \longrightarrow \ $   $\left(  M \rtimes L,\alpha_{\rtimes} \right)$, determines a homomorphism of Hom-Lie algebras  $f=\pi \circ h:\left(
K,\alpha_{K}\right)  \longrightarrow\left(  L,\alpha_{L}\right)$ and a derivation
$d=\theta \circ h:\left(  K,\alpha_{K}\right)  \longrightarrow\left(  M,\alpha
_{M}\right)$.
\end{proposition}
{\it Proof.}
Define $h:(K,\alpha_K)\to (M \rtimes L,\alpha_{\rtimes})$ by $h(x) = \left(  d\left(  x\right) ,f\left(  x\right)  \right)$, $x \in K$. Then everything can be readily checked. \rdg
\medskip

By taking $(K,\alpha_K)=(L,\alpha_L)$ and $f=\id_L$ we get:
\begin{corollary} Let $(M,\alpha_M)$ be a Hom-module over a Hom-Lie algebra $(L,\alpha_L)$.
The set of derivations from $(L,\alpha_L)$ to $(M,\alpha_M)$ is in a bijective correspondence with the set of
 homomorphisms $h : (L,\alpha_L) \to (M \rtimes L,\alpha_{\rtimes})$ such that  $\pi \circ h = \id_L$.
\end{corollary}

\begin{theorem} \label{sucesion}
Let $0 \to (N,\alpha_N) \stackrel{i}\to (K,\alpha_K)
\stackrel{\pi} \to (L,\alpha_L) \to 0$ be a short exact sequence of Hom-Lie algebras and  $(M,\alpha_M)$ a Hom-module over $(L,\alpha_L)$
 (and so a Hom-module over $(K,\alpha_K)$ via $\pi$ )
 such that  the Hom-action satisfies the condition ${^{\alpha_L(l)}} m = {^l} m$,  $l \in L$ and $m \in M$ (e.g. see Example \ref{ejemplo 2} {\it c)}). Denote by $(N^{\rm ab},\overline{\alpha}_{N})$ the quotient of $(N,\alpha_N)$ by the ideal $([N,N],\alpha_{N\mid})$.   Then $(N^{\rm ab},\overline{\alpha}_{N})$ has a Hom-module structure over $(L,\alpha_L)$ and  there is a natural exact sequence of vector spaces
\[
0 \to \Der_{\alpha}(L, M) \stackrel{\Delta}\longrightarrow \Der_{\alpha}(K, M) \stackrel{\rho}\longrightarrow \hom_L(N^{\rm ab},M),
\]
 where $\hom_L(N^{\rm ab},M)=\{ f : (N^{\rm ab},\overline{\alpha}_{N}) \to (M, \alpha_M) \mid f({^l}\overline{n}) = {^l}f(\overline{n}) \}$.
\end{theorem}
{\it Proof.} It is easy to see that the equality ${}^l\overline{n}=i^{-1}[x_l, i(n)]$, where $n\in N$, $l\in L$ and $x_l\in K$ such that $\pi(x_l)=l$, defines a Hom-module structure over $(L,\alpha_L)$ on $(N^{\rm ab},\overline{\alpha}_{N})$.

   Let $\Delta(d) = d \circ \pi$, for $d \in \Der_{\alpha}(L,M)$. Obviously $\Delta$ is injective. For any $\delta \in \Der_{\alpha}(K,M)$,  $\delta \circ i : (N,\alpha_N) \to (M, \alpha_M)$ is a homomorphism of Hom-Lie algebras that vanishes on $[N,N]$ and so induces a homomorphism of abelian Hom-Lie algebras $\rho : (N^{\rm ab}, \overline{\alpha}_N) \to (M, \alpha_M)$. Now the remaining details are straightforward. \rdg
\medskip

Let us note that, of course the results above recover the well-known classical facts on semi-direct product of Lie algebras (see e.g. \cite{HS}).

\subsection{Homology} \

\noindent The homology of Hom-Lie algebras, generalizing the classical Chevalley-Eilenberg homology of Lie algebras, is constructed in  \cite{Yau1, Yau2} (see also \cite{CaInPa}). Let us recall that the homology $ H_{*}^{\alpha}\left(  L,M\right)$ of a Hom-Lie algebra $(L,  \alpha_{L})$ with coefficients in a Hom-module $(M,\alpha_{M})$ over $(L,  \alpha_{L})$ is defined as the homology of the chain complex $(C_{*}^{\alpha}(L,M), d_*)$, where
\[
C_{n}^{\alpha}(L,M)=M\otimes \wedge^n L,\quad n\geq 0
\]
and the boundary map $d_{n}:C_{n}^{\alpha}(L,M)  \longrightarrow C_{n-1}^{\alpha}(L,M)$, $n\geq 1$, is given by
\begin{align*}
d_{n}&(  m\otimes x_{1}\wedge\cdots\wedge
x_{n})  =\overset{n}{\underset{i=1}
{\displaystyle\sum}}
(-1)^{i} \ \ {}^{x_{i}} m \otimes\alpha_{L}(x_{1})
\wedge\!\cdots\!\wedge\widehat{\alpha_{L}(x_{i})}\wedge\!\cdots\!\wedge \alpha_{L}(x_{n}) +\\
&\!\!\underset{1\leq i<j\leq n}{\sum}\!\!\left(  -1\right)
^{i+j}\alpha_{M}\left(  m\right)  \otimes\left[  x_{i},x_{j}\right]
\wedge\alpha_{L}\left(  x_{1}\right)  \wedge\!\cdots\!\wedge\widehat{\alpha
_{L}\left(  x_{i}\right)  }\wedge\!\cdots\!\wedge\widehat{\alpha_{L}\left(
x_{j}\right)  }\wedge\!\cdots\!\wedge\alpha_{L}\left(  x_{n}\right).
\end{align*}
As usual $\widehat{\alpha_{L}(x_{i})}$ means that the variable $\alpha_{L}(x_{i})$ is omitted.

Let us remark that for a Lie algebra $L$ and an $L$-module $M$, the chain complex
$C_{*}^{\alpha}(L,M)$ is exactly the Chevalley-Eilenberg complex that defines the
Lie algebra homology of $L$ with coefficients in the $L$-module $M$.

Easy computations of low-dimensional  cycles and boundaries provide the following results:
\[
H_{0}^{\alpha
}\left(  L,M\right)  ={\Ker\left(  d_{0}\right) }/{\Im\left(  d_{1}\right)  }={M}/{^{L}M} \ ,
\]
where $^{L}M=\left\{  {}^x m \mid m\in M, \ x\in L\right\}$. Moreover, if $(M,\alpha_M)$ is a trivial Hom-module over $(L,\alpha_L)$, i.e.  ${}^x m=0$ for all $x\in L$ and $m\in M$, then
\[
H_{1}^{\alpha}\left(  L,M\right)={\Ker\left(  d_{1}\right)  }/{\Im
\left(  d_{2}\right)  }  ={(M\otimes L)}/{\big(\alpha_M\left(  M\right)  \otimes\left[
L,L\right]\big)  }.
\]
In particular, if $M = \mathbb{K}$, then $H_{1}^{\alpha}\left(  L,\mathbb{K} \right)  ={L}/{\left[
L,L\right]  }$.

Below we use the notation $H_{n}^{\alpha}\left(  L \right) $ for $H_{n}^{\alpha}\left(  L,\mathbb{K} \right)$.

\section{Non-abelian tensor product of Hom-Lie algebras}

In this section we introduce a non-abelian tensor product of Hom-Lie algebras which generalizes the non-abelian tensor product of Lie algebras \cite{El2}, and study its properties.

\begin{definition}
Let $(M,\alpha_M)$ and $(N,\alpha_N)$ be Hom-Lie algebras with Hom-actions on each other. The Hom-actions are  said to be compatible if
\[
^{(^mn)} m'=[m',^nm] \quad \text{and} \quad ^{(^nm)}n'=[n',^mn]
\]
for all $m,m'\in M$ and $n,n'\in N$.
\end{definition}

\begin{example}
If  $(H,\alpha_H)$ and $(H',\alpha_{H'})$ both are ideals of a Hom-Lie algebra $(L,\alpha_L)$, then the Hom-actions of $(H,\alpha_H)$ and $(H',\alpha_{H'})$ on each other, considered in  Example \ref{ejemplo 2} {a)}, are compatible.
\end{example}

Let $(M,\alpha_M)$ and $(N,\alpha_N)$ be Hom-Lie algebras acting on each other compatibly.
Consider the Hom-vector space $(M\otimes N,\alpha_{M\otimes N})$  given by the tensor product $M\otimes N$ of the underlying vector spaces and the linear map $\alpha_{M \otimes  N} : M \otimes N \to M \otimes  N$, $\alpha_{M \otimes N}(m \otimes n) = \alpha_M(m) \otimes \alpha_N(n)$. Denote by $D(M,N)$ subspace of $M \otimes N$ generated by all elements of the form
\begin{enumerate}
\item[{\it a)}] $[m,m']\otimes \alpha_N(n) -\alpha_M(m)\otimes {}^{m'}n +\alpha_M(m') \otimes {}^mn$,
\item[{\it b)}] $\alpha_M(m)\otimes [n,n'] - {}^{n'}m\otimes\alpha_N(n)+{}^{n}m\otimes\alpha_N(n')$,
\item[{\it c)}] ${}^{n}m\otimes {}^{m}n$,
\item[{\it d)}]${}^{n}m\otimes {}^{m'}n' +{}^{n'}m'\otimes {}^{m}n$,
\item[{\it e)}]$[{}^{n}m,{}^{n'}m']\otimes \alpha_N({}^{m''}n'')+[{}^{n'}m',{}^{n''}m'']\otimes \alpha_N({}^{m}n)
+[{}^{n''}m'',{}^{n}m]\otimes \alpha_N({}^{m'}n')$,
\end{enumerate}
for $m,m',m''\in M$ and $n,n',n''\in N$.

\begin{proposition}
The quotient vector space $(M\otimes N)/D(M,N)$ with the product
\begin{equation}\label{pr_in_tensor}
[m\otimes n, m'\otimes n']=-{}^{n}m\otimes {}^{m'}n'
\end{equation}
and together with the endomorphism $(M\otimes N)/D(M,N)\to (M\otimes N)/D(M,N)$ induced by $\alpha_{M\otimes N}$, is a Hom-Lie algebra.
\end{proposition}
{\it Proof.}
It is clear that $\alpha_{M\otimes N}$ preserves the elements of $D(M,N)$ as well as the product defined by (\ref{pr_in_tensor}).
Routine calculations show that this product is compatible with the defining relations of $(M\otimes N)/D(M,N)$ and can be extended from generators to any elements. Since the actiones of $(M,\alpha_M)$ and $(N,\alpha_N)$ on each other are compatible, it follows by direct calculations that the product (\ref{pr_in_tensor}) satisfies the skew-symmetry and the Hom-Jacobi identity. \rdg

\begin{definition}
The above  Hom-Lie algebra structure on $(M\otimes N)/D(M,N)$ is called the non-abelian tensor product of Hom-Lie algebras  $(M,\alpha_M)$ and $(N,\alpha_N)$ (or Hom-Lie tensor product for short). It will be denoted by $(M\star N, \alpha_{M\star N})$ and the equivalence class of $m\otimes n$ will be denoted by $m\star n$.
\end{definition}

\begin{remark}
Note that if $\alpha_M=\id_M$ and $\alpha_N=\id_N$ then $M\star N$  is the non-abelian tensor product of Lie algebras $M$ and $N$ given in \cite{El2} (see also \cite{El1, InKhLa}).
\end{remark}

The Hom-Lie tensor product can also be defined by a universal property in the following way.

\begin{definition}\label{paring} Let $(M,\alpha_M)$ and $(N,\alpha_N)$ be Hom-Lie algebras acting on each other.
For any Hom-Lie algebra $(L,\alpha_L)$, a bilinear map $h : (M\times N,\alpha_M\times\alpha_N) \to (L,\alpha_L)$ is said to be a Hom-Lie pairing if the following properties are satisfied:
\begin{enumerate}
\item [a)] $h([m,m'],\alpha_N(n)) = h(\alpha_M(m), {}^{m'} n) - h(\alpha_M(m'), {}^m n)$,
 \item [b)] $h(\alpha_M(m),[n,n']) = h({}^{n'} m, \alpha_N(n))  - h({}^n m, \alpha_N(n'))$,
 \item [c)] $h({}^n  m, {}^{m'} n') = - [h(m,n), h(m',n')]$,
 \item [d)] $h  \circ (\alpha_M \times \alpha_N) = \alpha_L  \circ h$,
\end{enumerate}
for all $m, m' \in M$, $n, n' \in N$.
\end{definition}

\begin{example}\
\begin{enumerate}
\item[a)] If $\alpha_L=\id_L, \alpha_M= \id_M$ and $\alpha_N= \id_N$, then Definition \ref{paring} recovers the definition of Lie paring given in \cite{El1}.

\item[b)] Let $(M,\alpha_M)$ and $(N,\alpha_N)$ be ideals of a Hom-Lie algebra $(L,\alpha_L)$, then the bilinear map  $h : (M\times N,\alpha_M\times\alpha_N) \to (M \cap N,\alpha_{M \cap N})$, given by $h(m,n)=[m,n]$, is a Hom-Lie pairing.
    \end{enumerate}
\end{example}

\begin{definition}
A Hom-Lie pairing  $h : (M\times N,\alpha_M\times\alpha_N) \to (L,\alpha_L)$ is said to be universal if for any other  Hom-Lie pairing $h' : (M\times N,\alpha_M\times\alpha_N) \to (L',\alpha_{L'})$ there is a unique  homomorphism of Hom-Lie algebras $\theta : (L,\alpha_L) \to (L',\alpha_{L'})$ such that
$\theta  \circ h = h'$.
\end{definition}

Clearly, if $h$ is universal, then $(L,\alpha_{L})$ is determined up to isomorphism by $(M,\alpha_M), (N,\alpha_N)$ and the Hom-actions. Moreover, it is straightforward to show the following

\begin{proposition} Let $(M,\alpha_M)$ and $(N,\alpha_N)$ be Hom-Lie algebras acting on each other compatibly. The map
\[
h : (M\times N,\alpha_M\times\alpha_N) \to (M\star N, \alpha_{M\star N}),\quad (m, n) \mapsto m \star n
\]
is a universal Hom-Lie paring.
\end{proposition}

The Hom-Lie tensor product is symmetric in the sense of the following isomorphism of Hom-Lie algebras
\[
(M\star N, \alpha_{M\star N})\overset{\approx}{\longrightarrow} (N\star M, \alpha_{N\star M}), \quad m\star n\mapsto n\star m.
\]
This follows by the fact that $h : (M\times N,\alpha_M\times\alpha_N) \to (N\star M, \alpha_{N\star M})$, $(m, n)\mapsto n\star m$ is a Hom-Lie pairing and the universal property of $(M\star N, \alpha_{M\star N})$ thus yields a homomorphism $(M\star N, \alpha_{M\star N})\to (N\star M, \alpha_{N\star M})$, the inverse of which is defined similarly.

\

Sometimes the Hom-Lie tensor product can be described as the tensor product of vector spaces. In particular, we have the following

\begin{proposition}\label{proposition 5.10}
If the Hom-Lie algebras $(M,\alpha_M)$ and $(N, \alpha_N)$ act trivially on each other and both $\alpha_M$, $\alpha_N$ are epimorphisms, then there is an isomorphism of abelian Hom-Lie algebras
\[
(M\star N, \alpha_{M\star N})\approx (M^{ab}\otimes N^{ab}, \alpha_{M^{ab}\otimes N^{ab}}),
\]
where $M^{ab}=M/[M,M]$, $N^{ab}=N/[N,N]$ and $\alpha_{M^{ab}\otimes N^{ab}}$ is induced by $\alpha_M$ and $\alpha_N$.
\end{proposition}
{\it Proof.}
Since the Hom-actions are trivial, the relation (\ref{pr_in_tensor}) enables us to see that $(M\star N, \alpha_{M\star N})$ is an abelian Hom-Lie algebra. Further, since $\alpha_M$ and $\alpha_N$ are epimorphisms, the defining relations of the Hom-Lie tensor product say that the vector space $M\star N$ is the quotient of $M\otimes N$ by the relations
$[m,m']\otimes n=0=m\otimes [n,n']$ for all $m,m'\in M$,  $n,n'\in N$. The later is isomorphic to $M^{ab}\otimes N^{ab}$ and this isomorphism commutes with the endomorphisms $\alpha_{M^{ab}\otimes N^{ab}}$ and $\alpha_{M\star N}$. \rdg
\medskip

The Hom-Lie tensor product is functorial in the following sense: if $f:(M,\alpha_M)$ $\to (M',\alpha_{M'})$ and $g:(N,\alpha_N)\to (N',\alpha_{N'})$ are homomorphisms of Hom-Lie algebras together with compatible Hom-actions of $(M, \alpha_M)$ (resp. $(M', \alpha_{M'})$) and $(N, \alpha_N)$ (resp. $(N', \alpha_{N'})$) on each other such that $f$, $g$ preserve these Hom-actions, that is
\[
f({}^nm)={}^{g(n)}f(m), \quad g({}^mn)={}^{f(m)}g(n), \qquad m\in M,\ n\in N,
\]
then there is a homomorphism of Hom-Lie algebras
\[
f\star g:(M\star N,\alpha_{M\star N})\to (M'\star N',\alpha_{M'\star N'})
\]
 defined by $(f\star g)(m\star n)=f(m)\star g(n)$.

\begin{proposition}\label{exact-tensor-1}
Let $0\to (M_1,\alpha_{M_1})\overset{f}{\to}(M_2,\alpha_{M_2})\overset{g}{\to}(M_3,\alpha_{M_3})\to 0$ be a short exact sequence of Hom-Lie algebras. Let $(N,\alpha_{N})$ be a Hom-Lie algebra together with compatible Hom-actions of $(N,\alpha_{N})$ and $(M_i,\alpha_{M_i})$ $(i=1,2,3)$ on each other and $f$, $g$ preserve these Hom-actions.
Then there is an exact sequence of Hom-Lie algebras
\[
(M_1\star N,\alpha_{M_1\star N})\overset{f\star \id_N}{\longrightarrow}(M_2\star N,\alpha_{M_2\star N})\overset{g\star \id_N}{\longrightarrow}(M_3\star N,\alpha_{M_3\star N})\longrightarrow 0.
\]
\end{proposition}
{\it Proof.}
Clearly $g\star \id_N$ is an epimorphism and $\Im(f\star \id_N) \subseteq \Ker(g\star \id_N)$. Now $\Im(f\star \id_N)$ is generated by all elements of the form $f(m_1)\star n_1$ with $m_1\in M_1$, $n_1\in N$ and it is an ideal in $(M_2\star N,\alpha_{M_2\star N})$ since we have
\[
[f(m_1)\star n_1, m_2\star n_2]=-f({}^{n_1}m_1)\star {}^{m_2}n_2\in \Im(f\star \id_N)
\]
for any generator $m_2\star n_2\in M_2\star N$. Thus, $g\star \id_N$ yields a factorization
\[
\xi:\big((M_2\star N) /\Im(f\star \id_N), \overline{\alpha}_{M_2\star N}\big)\to (M_3\star N,\alpha_{M_3\star N}).
\]
In fact this is an isomorphism of Hom-Lie algebras with the inverse map
\[
\xi':(M_3\star N,\alpha_{M_3\star N}) \to \big((M_2\star N) /\Im(f\star \id_N), \overline{\alpha}_{M_2\star N}\big)
\]
given on generators by $\xi'(m_3\star n)=\overline{m_2\star n}$, where $m_2\in M_2$ such that $g(m_2)=m_3$.
The remaining details are straightforward calculations and we leave to the reader. \rdg

\begin{proposition}\label{exact-tensor-2}
 If $(M,\alpha_M)$ is an ideal of a Hom-Lie algebra $(L, \alpha_L)$, then there is an exact sequence of Hom-Lie
algebras
\[
\big((M\star L)\rtimes (L\star M), \alpha_{\rtimes}\big) \overset{\sigma} {\longrightarrow}(L\star L,\alpha_{L\star L})\overset{\tau}{\longrightarrow}({L}/{M}\star {L}/{M},\alpha_{{L}/{M}\star {L}/{M}})\longrightarrow 0.
\]
\end{proposition}
{\it Proof.}
First we note that $\tau$ is the functorial homomorphism induced by the projection $(L,\alpha_L) \twoheadrightarrow (L/M,\alpha_{L/M})$ and clearly it is surjective. Let $\sigma':(M\star L,\alpha_{M\star L})\to (L\star L,\alpha_{L\star L})$ and  $\sigma'':(L\star M,\alpha_{L\star M}) \to (L\star L,\alpha_{L\star L})$ be the functorial homomorphisms induced by the inclusion $(M,\alpha_M)\hookrightarrow (L,\alpha_L)$ and by the identity map $(L,\alpha_L)\to (L,\alpha_L)$. Let $\sigma (x,y)=\sigma'(x)+\alpha_{L \star M} \circ \sigma''(y)$ for all $x\in M\star L$ and $y\in L\star M$. It is straightforward to see that $\sigma$ is a homomorphism of  Hom-Lie algebras and $\tau \circ \sigma$ is the trivial homomorphism. Clearly $\Im(\sigma)$ is generated by the elements $m\star l$ and $\alpha_L(l)\star \alpha_M(m)$ for $m\in M$, $l\in L$ and, by the formula (\ref{pr_in_tensor}), it is an ideal of $(L\star L,\alpha_{L\star L})$. Let us define a homomorphism of Hom-Lie algebras $\tau': ({L}/{M}\star {L}/{M},\alpha_{{L}/{M}\star {L}/{M}})\to (L\star L,\alpha_{L\star L})/\Im(\sigma)$ by $\tau'(\overline{l}\star \overline{l'})=\overline{l\star l'}$, $l,l'\in L$. It is easy to see that $\tau'$ is well-defined and it has an inverse homomorphism induced by $\tau$. \rdg

\begin{lemma} \label{action-on-tensor}
Let  $(M,\alpha_M)$ and $(N,\alpha_N)$ be Hom-Lie algebras with compatible actions on each other.
 \begin{enumerate}
\item[a)] There are homomorphisms of Hom-Lie algebras

$\begin{array}{rl}
&\psi_M:(M\star N, \alpha_{M \star N}) \to (M, \alpha_{M }), \quad \psi_M(m\star n)= -{}^nm,\\
&\psi_N:(M\star N, \alpha_{M \star N}) \to (N, \alpha_N),\quad  \psi_N(m\star n)= {}^mn.
\end{array}$
\item[b)]
There is a Hom-action of  $(M, \alpha_{M })$ (resp. $(N, \alpha_{N })$) on  the Hom-Lie tensor product ($M\star N, \alpha_{M \star N}$) given, for all $m,m'\in M$, $n,n'\in N$, by
\[
\begin{array}{rl}
&{}^{m'}(m\star n)=[m',m]\star \alpha_N(n)+\alpha_M(m)\star {}^{m'}n \\
\big( \text{resp.} \ &{}^{n'}(m\star n)={}^{n'}m\star \alpha_N(n)+\alpha_M(m)\star [n',n] \big)
\end{array}
\]
\item[c)] $\Ker(\psi_M)$ (resp. $\Ker(\psi_N)$) is contained in the center of $(M\star N, \alpha_{M \star N})$.
\item[d)] The induced Hom-action of $\Im(\psi_1)$ (resp. $\Im(\psi_2)$) on $\Ker(\psi_1)$ (resp. $\Ker(\psi_2)$) is trivial.
\item[e)]  $\psi_M$ and $\psi_N$ satisfy the following properties for all $m, m' \in M$, $n, n' \in N$:
\begin{enumerate}
\item[i)] $\psi_M(^{m'}(m \star n)) = [\alpha_M(m'), \psi_M(m \star n)]$,
\item[ii)] $\psi_N(^{n'}(m \star n)) = [\alpha_N(n'), \psi_N(m \star n)]$,
\item[iii)] ${^{\psi_M(m \star n)}}(m' \star n') = [\alpha_{M \star N}(m \star n), m' \star n'] = {^{\psi_N(m \star n)}} (m' \star n')$.
\end{enumerate}
\end{enumerate}
\end{lemma}
{\it Proof.}
Everything can be readily checked thanks to the compatibility conditions and the relation (\ref{pr_in_tensor}). \rdg

\begin{remark}
If $\alpha_M = id_M$ and $\alpha_N = id_N$, then both $\psi_M$ and $\psi_N$ in Theorem \ref{action-on-tensor} are crossed modules of Lie algebras (see \cite{El1}).
\end{remark}

\begin{definition}
A Hom-Lie algebra  $(L, \alpha_L)$ is said to be perfect if $L=[L, L]$.
\end{definition}

\begin{theorem} \label{sucesion exacta}
Let $(M,\alpha_M)$ be an ideal of a perfect Hom-Lie algebra $(L,\alpha_L)$. Then there is an exact sequence of vector spaces
 \[
 \Ker(M\star L\overset{\psi_{M}\ }\longrightarrow M)\to H_2^{\alpha}(L)\to H_2^{\alpha}({L}/{M})\to {M}/{[L,M]}\to  0
 \]
\end{theorem}
{\it Proof.} Thanks to Proposition \ref{exact-tensor-2} there is a commutative diagram of Hom-Lie algebras with exact rows
{\footnotesize\[
 \xymatrix{
  & \big((M\star L) \rtimes (L\star M), \alpha_{\rtimes}\big) \ar[r]  \ar[d]^{\psi}& (L \star L,\alpha_{L \star L}) \ar[r]^{\pi \star \pi \ \ \ \ \ \ \ \ } \ar[d]^{\psi_L} & ({L}/{M} \star {L}/{M}, {\alpha}_{{L}/{M} \star {L}/{M}}) \ar[r]  \ar[d]^{{\psi_{L/M}}}& 0\\
0 \ar[r] & (M, \alpha_M) \ar[r] & (L,\alpha_L) \ar[r]^{\pi} & ({L}/{M},{\alpha}_{L/M}) \ar[r] & 0,
}
\]}
where $\psi((m_1,l_1),(l_2,m_2))=[m_1,l_1]+[\alpha_L(l_2),\alpha_M(m_2)]$. Then, by using the Snake Lemma,  the assertion follows from Remark \ref{H2} below and the fact that there is a surjective
map $\Ker(\psi) \to \Ker(\psi_{M})$. \rdg

\begin{remark}
If $\alpha_L = \id_L$, then the exact sequence in Theorem \ref{sucesion exacta} is part of the six-term exact sequence in \cite{El2}.
\end{remark}

\section{Application in universal ($\alpha$-)central extensions of Hom-Lie algebras} \label{section}

In this section we complement by new results the investigation of universal central extensions of Hom-Lie algebras done in \cite{CaInPa}. We also describe universal ($\alpha$-)central extensions via Hom-Lie tensor product.

\begin{definition} \label{alfacentral} A central (resp. $\alpha$-central) extension of a Hom-Lie algebra $(L, \alpha_L)$ is an exact sequence of Hom-Lie algebras
 \[
 (\mathfrak{K}): \ \ 0 \longrightarrow (M, \alpha_M)  \longrightarrow (K,\alpha_K) \stackrel{\pi} \longrightarrow (L, \alpha_L) \longrightarrow 0
 \]
such that $[M, K] = 0 $, i.e. $M \subseteq Z(K)$ (resp. $[\alpha_M(M), K] = 0$, i.e. $\alpha_M(M) \subseteq Z(K)$).

A central extension $(\mathfrak{K})$ is called universal central (resp. universal $\alpha$-central) extension if, for every central (resp. $\alpha$-central) extension $(\mathfrak{K'})$ of $(L,\alpha_L)$ there exists one and only one homomorphism of Hom-Lie algebras $h : (K,\alpha_K) \to (K',\alpha_{K'})$ such that $\pi' \circ h = \pi$.
\end{definition}

\begin{remark}
Obviously every  central extension is an $\alpha$-central extension and these notions coincide when $\alpha_M = \id_M$. On the other hand, every universal $\alpha$-central extension  is a  universal  central extension and these notions coincide when $\alpha_M = \id_M$. Let us also observe that if a universal ($\alpha$-)central  extension exists then it is unique up to isomorphism.
\end{remark}

The category ${\tt HomLie}$ is an example of a semi-abelian category which does not satisfy universal central extension condition in the sense of \cite{CaTim}, that is, the composition of central extensions of Hom-Lie algebras is not central in general, but it is an $\alpha$-central extension (see Theorem \ref{teorema} {\it a)} below). This fact does not allow complete generalization of classical results to Hom-Lie algebras and the well-known properties of universal central extensions are divided between universal central and universal $\alpha$-central extensions of Hom-Lie algebras. In particular, the assertions in the following theorem are proved in \cite{CaInPa}.

\begin{theorem}\label{teorema} \
\begin{enumerate}
\item[a)] Let $(K,\alpha_K) \stackrel{\pi} \twoheadrightarrow (L, \alpha_L)$  and  $(F,\alpha_F) \stackrel{\rho} \twoheadrightarrow (K, \alpha_K)$ be  central extensions with  $(K, \alpha_K)$ a perfect Hom-Lie algebra. Then the composition extension $(F,\alpha_F)$ $\stackrel{\pi \circ \rho} \twoheadrightarrow (L, \alpha_L)$ is an $\alpha$-central extension.

    \item[b)] Let  $(K,\alpha_K) \stackrel{\pi} \twoheadrightarrow (L, \alpha_L)$ and $(K',\alpha_{K'}) \stackrel{\pi'} \twoheadrightarrow (L, \alpha_L)$ be two  central extensions of $(L,\alpha_L)$. If $(K,\alpha_K)$ is perfect, then there exists at most one homomorphism of Hom-Lie algebras  $f : (K,\alpha_K) \to (K', \alpha_{K'})$ such that $\pi' \circ f = \pi$.

\item[c)]  If $ (K,\alpha_K) \stackrel{\pi} \twoheadrightarrow (L, \alpha_L)$ is a universal $\alpha$-central extension, then  $(K,\alpha_K)$ is a perfect   Hom-Lie algebra and every central extension of $(K,\alpha_K)$ splits.

    \item[d)] If $(K,\alpha_K)$ is a perfect  Hom-Lie algebra and every  central extension of $(K,\alpha_K)$ splits, then any central extension $(K,\alpha_K) \stackrel{\pi} \twoheadrightarrow (L, \alpha_L)$ is a universal central extension.

\item[e)] A Hom-Lie algebra $(L, \alpha_L)$ admits a universal central extension if and only if $(L, \alpha_L)$ is perfect.
Furthermore, the kernel of the universal central extension is canonically isomorphic to the second homology $H_2^{\alpha}(L)$.

\item[f)] If $(K,\alpha_K) \stackrel{\pi} \twoheadrightarrow (L, \alpha_L)$ is a universal $\alpha$-central extension, then $H_1^{\alpha}(K) = H_2^{\alpha}(K) = 0$.

   \item[g)]   If $H_1^{\alpha}(K) = H_2^{\alpha}(K) = 0$, then any central extension  $(K,\alpha_K) \stackrel{\pi} \twoheadrightarrow (L, \alpha_L)$  is a universal  central extension.

\end{enumerate}
\end{theorem}

It follows from Lemma \ref{action-on-tensor} that for any Hom-Lie algebra $(L,\alpha_L)$ the homomorphism
 \[
 \psi:(L\star L, \alpha_{L \star L}) \twoheadrightarrow ([L,L], \alpha_{L \mid}), \quad  \psi(l \star l')=[l,l'],
 \]
 is a central extension of the Hom-Lie algebra $([L,L], \alpha_{L \mid})$.

\begin{theorem} \label{teor}
If $(L,\alpha_L)$ is a perfect Hom-Lie algebra,  then the central extension \linebreak $(L\star L, \alpha_{L \star L})\overset{\psi}\twoheadrightarrow (L, \alpha_L)$ is the universal central extension of $(L,\alpha_L)$.
\end{theorem}
{\it Proof.} Let $(C, \alpha_C) \overset{\phi}\twoheadrightarrow (L, \alpha_L)$ be a central extension of $(L, \alpha_L)$. Since $\Ker(\phi)$ is in the center of $(C,\alpha_C)$, we get a well-defined homomorphism of Hom-Lie algebras $f:(L\star L, \alpha_{L\star L})\to (C, \alpha_C)$ given on generators by $f(l\star l')=[c_l,c_{l'}]$, where $c_l$ and $c_{l'}$ are any elements in $\phi^{-1}(l)$ and $\phi^{-1}(l')$, respectively. Obviously $\phi \circ f= \psi$ and $f \circ \alpha_{L \star L} = \alpha_C \circ f$, having in mind that $\alpha_C(c_l) \in \phi^{-1}(\alpha_L(l))$ for all $l \in L$.  Since $L$ is perfect, then by equality (\ref{pr_in_tensor}), so is $L\star L$. Hence the homomorphism $f$ is unique by Theorem \ref{teorema} {\it b)}. \rdg

\begin{remark} \label{H2}
If the Hom-Lie algebra $(L,\alpha_L)$ if perfect, by Theorem \ref{teorema} {e)} we have that $H_2^{\alpha}(L)\approx \Ker(L\star L\overset{\psi}\to L )$.
\end{remark}

Now we obtain a condition for the existence of the universal $\alpha$-central extensions. We need the following notion.

\begin{definition}
A Hom-Lie algebra  $(L, \alpha_L)$ is said to be $\alpha$-perfect if $L=[\alpha_L(L),$ $\alpha_L(L)]$.
\end{definition}

\begin{example}
Consider the situation when  the ground field $\K$ is the field of complex numbers. Let $L$ be the three-dimensional vector space with basis $\{e_1,$ $e_2,e_3\}$. Define product in $L$ by $[e_1,e_2]= - [e_2,e_1] = e_3$, $[e_2,e_3]= - [e_3,e_2] = e_1$, $[e_3,e_1]= - [e_1,e_3] = e_2$ and zero elsewhere. Take the endomorphism $\alpha_L:L\to L$ represented by the matrix $\left( \begin{array}{ccc} \frac{\sqrt{2}}{2} & 0 & \frac{\sqrt{2}}{2} \\ 0 & -1 & 0 \\ \frac{\sqrt{2}}{2} & 0 & - \frac{\sqrt{2}}{2} \end{array} \right)$. Then $(L,\alpha_L)$ is  an $\alpha$-perfect Hom-Lie algebra.
\end{example}

 \begin{remark}\label{alfa perfecta} \
 \begin{enumerate}
 \item[a)] When $\alpha_L= \id_L$, the notions of perfect and $\alpha$-perfect Hom-Lie algebras are the same.

 \item[b)] Obviously, if $\left(  L,\alpha_{L}\right)$ is an $\alpha$-perfect Hom-Lie algebra, then it is perfect. Nevertheless the converse is not true in general. For example, the three-dimensional (as a vector space) Hom-Lie  algebra $\left(  L, \alpha_{L}\right)$  with linear basis $\{e_1, e_2, e_3\}$, product given by $[e_1,e_2]=-[e_2,e_1]=e_3$, $[e_1,e_3]=-[e_3,e_1]=e_2$, $[e_2,e_3]=-[e_3,e_2]=e_1$ and zero elsewhere, and endomorphism $\alpha_L=0$ is perfect, but it is not $\alpha$-perfect.

     \item[c)] If  $\left(  L,\alpha_{L}\right)$  is $\alpha$-perfect, then $L = \alpha_L(L)$, i.e. $\alpha_L$ is surjective. Nevertheless the converse is not true. For instance, consider the two-dimensional (as a vector space) Hom-Lie algebra with linear basis $\{e_1, e_2\}$, bracket given by $[e_1, e_2] = - [e_2, e_1] = e_2$ and zero elsewhere, and endomorphism $\alpha_L$ represented by the matrix $\left( \begin{array}{cc} 1 & 0 \\ 0 & 2 \end{array} \right)$. Obviously the endomorphism $\alpha_L$ is surjective, but $[\alpha_L(L), \alpha_L(L)] = \langle \{e_2 \} \rangle$.

 \end{enumerate}
 \end{remark}

 \begin{lemma} \label{lema 10}
Let $(M, \alpha_M)  \rightarrowtail (K,\alpha_K) \stackrel{\pi} \twoheadrightarrow (L, \alpha_L)$ be a central extension and $(K,\alpha_K)$ be an $\alpha$-perfect Hom-Lie  algebra. Let $ (M',\alpha_{M'})  \rightarrowtail (K',\alpha_{K'}) \stackrel{\pi'} \twoheadrightarrow (L, \alpha_L)$ be an $\alpha$-central extension. Then there exists at most one homomorphism of Hom-Lie algebras  $f : (K,\alpha_K) \to (K', \alpha_{K'})$ such that $\pi' \circ f = \pi$.
\end{lemma}
{\it Proof.} Let us assume that there are homomorphisms  $f_1$ and $f_2$ such that $\pi' \circ f_1 = \pi = \pi' \circ f_2$. Then for any $k\in K$ we have $f_1(k) = f_2(k)+m'_k$,  for some $m'_k \in M'$. By using the condition $\alpha_{M'}(M') \subseteq Z(K')$ we have
 \begin{align*}
 f_1[\alpha_K(k_1),\alpha_K(k_2)]&=[\alpha_{K'}  f_1(k_1), \alpha_{K'}  f_1(k_2)] \\
 & = [\alpha_{K'} f_2(k_1) \!+\! \alpha_{K'}(m'_{k_1}), \alpha_{K'}  f_2(k_2)\! +\! \alpha_{K'}(m'_{k_2})]\\
 & = [\alpha_{K'} f_2(k_1) , \alpha_{K'} f_2(k_2)]\\
 & = f_2[\alpha_K(k_1),\alpha_K(k_2)]
 \end{align*}
 for any $k_1,k_2\in K$. This implies that $f_1=f_2$, since $(K,\alpha_K)$ is $\alpha$-perfect. \rdg

\begin{theorem}\label{alfa uce}
An $\alpha$-perfect Hom-Lie algebra admits a universal $\alpha$-central extension.
\end{theorem}
{\it Proof.}
Given an $\alpha$-perfect Hom-Lie algebra $(L,\alpha_L)$ we construct a universal $\alpha$-central extension
\begin{equation}\label{eq}
0 \to \Ker (u_{\alpha}) \to (\frak{uce}_{\alpha}(L), \widetilde{\alpha})  \stackrel{u_{\alpha}} \to (L, \alpha_L) \to 0
\end{equation}
as follows. We consider the quotient vector space  $\frak{uce}_{\alpha}(L) =\big({\alpha_L(L) \wedge \alpha_L(L)}\big)/{I_L}$, where $I_L$ is the vector subspace of $\alpha_L(L) \wedge \alpha_L(L)$ spanned by the elements of the form
$$-[x_1,x_2] \wedge \alpha_L(x_3) + [x_1,x_3] \wedge \alpha_L(x_2) -  [x_2,x_3] \wedge \alpha_L(x_1)$$
for all $x_1, x_2, x_3 \in L$. Here we observe that every summand of the form $[x_1,x_2] \wedge \alpha_L(x_3)$  is an element of $\alpha_L(L) \wedge \alpha_L(L)$, since  $L$  is $\alpha$-perfect and so $[x_1,x_2] \in L = [\alpha_L(L), \alpha_L(L)] \subseteq \alpha_L(L)$.
We denote by $\{\alpha_L(x_1), \alpha_L(x_2)\}$ the equivalence class of $\alpha_L(x_1) \wedge \alpha_L(x_2)$. The product in  $\frak{uce}_{\alpha}(L)$ is defined by
\[[\{\alpha_L(x_1),\alpha_L(x_2)\}, \{\alpha_L(y_1),\alpha_L(y_2)\}] = \{ [\alpha_L(x_1),\alpha_L(x_2)], [\alpha_L(y_1),\alpha_L(y_2)]\}
\]
and the endomorphism $\widetilde{\alpha} : \frak{uce}_{\alpha}(L) \to \frak{uce}_{\alpha}(L)$ is given by
\[
\widetilde{\alpha}(\{\alpha_L(x_1),\alpha_L(x_2)\}) = \{\alpha_L^2(x_1),\alpha_L^2(x_2)\}.
\]
The map $u_{\alpha}$  is defined by $u_{\alpha}(\{\alpha_L(x_1), \alpha_L(x_2)\})= [\alpha_L(x_1), \alpha_L(x_2)]$.
Straightforward calculations show that $(\frak{uce}_{\alpha}(L), \widetilde{\alpha})$ is indeed a Hom-Lie algebra and
$u_{\alpha}$ is a homomorphism of Hom-Lie algebras. Moreover, $u_{\alpha}$ is surjective, because  $(L,\alpha_L)$ is $\alpha$-perfect.

Obviously the sequence (\ref{eq}) is a central extension. Moreover, it is a universal $\alpha$-central extension.
Indeed, consider any $\alpha$-central extension $(M,
\alpha_M) \rightarrowtail (K, \alpha_K) \stackrel{\pi} \twoheadrightarrow (L,
\alpha_L)$. We define $\Phi : (\frak{uce}_{\alpha}(L), \widetilde{\alpha}) \to (K, \alpha_K)$ by
$\Phi(\{\alpha_L(x_1), \alpha_L(x_2)\})$ $= [\alpha_K(k_1),
\alpha_K(k_2)]$, where $k_1, k_2 \in K$ such that $\pi(k_1)=x_1$, $\pi(k_2)=x_2$. It is well defined  because of the equality $[\alpha_M(M) , K ]= 0$.
Moreover, direct calculations show that $\Phi$ is a homomorphism of Hom-Lie algebras and $\pi \circ \Phi=u_{\alpha}$. To prove the uniqueness of such $\Phi$, by Lemma \ref{lema 10} it is enough to check that $\frak{uce}_{\alpha}(L)$ is
$\alpha$-perfect. For this later we do the following calculations:
 \[
 [\widetilde{\alpha} \{\alpha_L(x_1), \alpha_L(x_2)\}, \widetilde{\alpha} \{\alpha_L(y_1), \alpha_L(y_2)\}] = \{[\alpha_L^2(x_1), \alpha_L^2(x_2)],[\alpha_L^2(y_1), \alpha_L^2(y_2)]\},
 \]
which implies that $[\widetilde{\alpha} (\frak{uce}_{\alpha}(L)),
\widetilde{\alpha}( \frak{uce}_{\alpha}(L))] \subseteq
\frak{uce}_{\alpha}(L)$. Conversely, having in mind that $L=[\alpha_L(L),\alpha_L(L)]$ and hence every element $x\in L$ can be written as $x=\displaystyle \sum_i \lambda_i [\alpha_L(l_{i_1}), \alpha_L(l_{i_2})] $ for some $\lambda_i\in \mathbb{K}$ and $l_{i_1}, l_{i_2}\in L$, we get
 \begin{align*}
 \{\alpha_L&(x_1), \alpha_L(x_2)\} = \left \{\alpha_L \left( \sum_i \lambda_i [\alpha_L(l_{i_1}), \alpha_L(l_{i_2})] \right), \alpha_L \left( \sum_j \lambda'_j [\alpha_L(l'_{j_1}), \alpha_L(l'_{j_2})] \right)\right \}\\
  &=\sum_{i,j} \lambda_i \lambda'_j \left\{ [\alpha^2_L(l_{i_1}),\alpha^2_L(l_{i_2})], [\alpha^2_L(l'_{j_1}),\alpha^2_L(l'_{j_2})] \right \}\\
  &= \sum_{i,j} \lambda_i \lambda'_j \left[ \{\alpha^2_L(l_{i_1}),\alpha^2_L(l_{i_2})\} , \{\alpha^2_L(l'_{j_1}),\alpha^2_L(l'_{j_2})\} \right] \\
  &=\sum_{i,j} \lambda_i \lambda'_j \left[ \widetilde{\alpha} \{\alpha_L(l_{i_1}),\alpha_L(l_{i_2})\}, \widetilde{\alpha} \{\alpha_L(l'_{j_1}),\alpha_L(l'_{j_2})\} \right] \in [\widetilde{\alpha} (\frak{uce}_{\alpha}(L)),\widetilde{\alpha} (\frak{uce}_{\alpha}(L))]
\end{align*}
for any $\{\alpha_L(x_1), \alpha_L(x_2)\}\in \frak{uce}_{\alpha}(L)$. \rdg

\begin{theorem}
If  $(L,\alpha_L)$ is an $\alpha$-perfect Hom-Lie algebra, then the homomorphism of
Hom-Lie algebras $\varphi:(\alpha_L(L)\star \alpha_L(L), \alpha_{\alpha_L(L) \star \alpha_L(L)})\twoheadrightarrow (L, \alpha_L)$ given by $\varphi(\alpha_L(l_1) \star \alpha_L(l_2))=[\alpha_L(l_1), \alpha_L(l_2)]$,  is the universal $\alpha$-central extension of $(L,\alpha_L)$.
Moreover, there is an isomorphism of Hom-Lie algebras
\[
(\alpha_L(L) \star \alpha_L(L), \alpha_{\alpha_L(L) \star \alpha_L(L)}) \approx (\frak{uce}_{\alpha}(L), \widetilde{\alpha}),\quad \alpha_L(l) \star \alpha_L(l') \mapsto \{\alpha_L(l),\alpha_L(l')\}.
\]
\end{theorem}
{\it Proof.} This is similar to the proof of Theorem \ref{teor} and we leave to the reader. \rdg

\section{Application in cyclic homology of Hom-associative algebras}\label{section 6}

Throughout this section we assume that $\K$ is a field of characteristic 0.

\begin{definition}
 By a Hom-associative algebra (see e.g. \cite{MS}) we mean a pair $(A,\alpha_A)$ consisting of a vector space $A$ and a linear map $\alpha_A:A\to A$, together with a linear map (multiplication) $A\otimes A\to A$, $a\otimes b\mapsto ab$,  such that
\begin{align*}
&\alpha_A(a)(bc) = (ab)\alpha_A(c),\\
& \alpha_A(ab)=\alpha_A(a)\alpha_A(b)
\end{align*}
for all $a,b,c\in A$.
\end{definition}

The Hom version of the classical cyclic bicomplex (see e.g. \cite{Lo})  is constructed in \cite{Yau1} and
the cyclic homology of a Hom-associative algebra is defined as the homology of its total complex. A reformulation of this cyclic homology via Connes's complex for Hom-associative algebra is also given in \cite[Proposition 4.7]{Yau1}. It follows that, given a Hom-associative algebra $(A,\alpha_A)$, the first cyclic homology $HC^{\alpha}_1(A)$ is the kernel of the homomorphism of vector spaces
\[
\psi:A\otimes A/J(A,\alpha)\to [A,A], \quad a\otimes b\mapsto ab-ba,
\]
where $[A,A]$ is the subspace of $A$ generated by the elements $ab-ba$, and $J(A,\alpha)$ is the subspace of $A\otimes A$ generated by the elements
\[
a\otimes b+b\otimes a \quad\text{and}\quad ab\otimes \alpha_A(c)- \alpha_A(a)\otimes bc+ca\otimes \alpha_A(b),
\]
for all $a,b,c\in A$

Given a Hom-associative algebra $(A,\alpha_A)$, then it is endowed with a Hom-Lie algebra structure with the induced product  $[a,b]=ab-ba$, $a,b\in A$ and the endomorphism $\alpha_A$. Moreover,
there is a Hom-Lie algebra structure on $A\otimes A/J(A,\alpha)$ given by
\[
[a\otimes b, a'\otimes b']=[a,b]\otimes [a',b']
\]
and the endomorphism induced by $\alpha_A$. We denote this Hom-Lie algebra by $(L^{\alpha}(A),\overline{\alpha}_A)$. In fact $(L^{\alpha}(A),\overline{\alpha}_A)$ is  the quotient of the Hom-Lie tensor product $(A\star A,\alpha_{A\star A})$ by the ideal generated by the elements
$a\star b+b\star a$ and $ab\star \alpha_A(c)- \alpha_A(a)\star bc+ca\star \alpha_A(b)$, for all $a,b,c\in A$.

\begin{definition} \label{alfa condition}
We say that a Hom-associative algebra $(A,\alpha_A)$ satisfies the $\alpha$-identity condition if
\begin{equation}\label{condition}
[A, \Im(\alpha_A-\id_A)]=0,
\end{equation}
where $[A, \Im(\alpha_A-\id_A)]$ is the subspace of $A$ spanned by all elements $ab-ba$ with $a\in A$ and $b\in \Im(\alpha_A-\id_A)$.
\end{definition}

Note that $\alpha$-identity condition is equivalent to the condition $[a,b]=[\alpha_A(a),b]$ for all $a,b\in A$.
\newpage

\begin{example} \
\begin{enumerate}
\item[a)] Any Hom-associative algebra $(A, \alpha_A)$ with $\alpha_A=\id_A$ (i.e. an associative algebra) satisfies  $\alpha$-identity condition.
\item[b)] Any commutative Hom-associative algebra $(A, \alpha_A)$ (i.e. $ab=ba$ for all $a,b\in A$) with $\alpha_A=0$ satisfies  $\alpha$-identity condition.
\item[c)] Consider the Hom-associative algebra $(A,\alpha_A)$, where as vector space $A$ is 2-dimensional with basis $\{e_1,e_2\}$,  the multiplication is given by $e_1e_1=e_2$ and zero elsewhere, $\alpha_A$ is represented by the matrix $\left( \begin{array}{cc} 1 & 0 \\ 1 & 1 \end{array} \right)$. Then
$(A,\alpha_A)$ satisfies  $\alpha$-identity condition.
\item[d)] Consider the Hom-associative algebra $(A,\alpha_A)$, where as vector space $A$ is 3-dimensional with basis $\{e_1,e_2,e_3\}$,  the multiplication is given by $e_1e_1=e_2$, $e_1e_2=e_3$, $e_2e_1=e_3$ and zero elsewhere, $\alpha_A$ is represented by the matrix $\left( \begin{array}{ccc} 1 & 0 &0 \\ 1 & 0 & 0 \\ 1 & 0 & 0 \end{array} \right)$. Then $(A,\alpha_A)$ satisfies  $\alpha$-identity condition.
\end{enumerate}
\end{example}

\begin{lemma}\label{lema 6.4}
Let $(A,\alpha_{ A})$  be a Hom-associative algebra.
\begin{enumerate}
\item[a)] There are Hom-actions of Hom-Lie algebras $(A,\alpha_A)$ and $(L^{\alpha}(A),\overline{\alpha}_A)$ on each other. Moreover, these Hom-actions are compatible if $(A,\alpha_A)$ satisfies  the $\alpha$-identity condition (\ref{condition}).
\item[b)] There is a short exact sequence of Hom-Lie algebras
\[
0\longrightarrow (HC^{\alpha}_1(A),\alpha_{HC})\overset{i}\longrightarrow (L^{\alpha}(A),\overline{\alpha}_A) \overset{\psi}\longrightarrow \big([A,A], \alpha_{A\mid}\big)\longrightarrow 0,
\]
where $(HC^{\alpha}_1(A),\alpha_{HC})$ is an abelian Hom-Lie algebra with $\alpha_{HC}$ induced by $\alpha_A$, $\alpha_{A\mid}$ is the restriction of $\alpha_A$ and $\psi(a\otimes b)=[a,b]$.

\item[c)] The induced Hom-action of $(A,\alpha_A)$ on  $(HC^{\alpha}_1(A),\alpha_{HC})$ is trivial.  Moreover, if $(A,\alpha_A)$ satisfies the $\alpha$-identity condition (\ref{condition}), then both $i$ and $\psi$ preserve the Hom-actions of the Hom-Lie algebra $(A,\alpha_A)$.
\end{enumerate}
\end{lemma}
{\it Proof.}
{\it a)} The Hom-action of $(A,\alpha_A)$ on $(L^{\alpha}(A),\overline{\alpha}_A)$ is given by
\[
{}^{a'}(a\otimes b) = [a',a]\otimes \alpha_A( b) + \alpha_A(a)\otimes [a',b],
\]
while  the Hom-action of $(L^{\alpha}(A),\overline{\alpha}_A)$ on $(A,\alpha_A)$ is defined by
\[
{}^{(a\otimes b)}{a'}=[[a,b],a']
\]
for all $a',a,b\in A$. Straightforward calculations show that these are indeed Hom-actions of Hom-Lie algebras, which are compatible if  $(A,\alpha_A)$ satisfies  $\alpha$-identity condition (\ref{condition}).

{\it b)} and {\it c)} are immediate consequences of the definitions above. \rdg

By complete analogy to the Dennis-Stein generators \cite{DeSt}, we define the first Milnor cyclic homology for Hom-associative algebras as follows.

\begin{definition}
Let $(A,\alpha_A)$ be a Hom-associative algebra. The first Milnor cyclic homology $HC_1^M(A,\alpha_A)$ is the quotient vector space of $A\otimes A$ by the relations
\begin{align*}
& a\otimes b +b\otimes a=0,\\
& ab\otimes \alpha_A(c)- \alpha_A(a)\otimes bc+ca\otimes \alpha_A(b)=0,\\
& \alpha_A(a)\otimes bc -\alpha_A(a)\otimes cb =0
\end{align*}
for all $a,b,c\in A$.
\end{definition}

Of course for $\alpha_A=\id_A$ this is the definition of the first Milnor cyclic homology of the associative algebra
$A$ in the sense of \cite{Lo} (see also \cite{InKhLa}).
Note also that $HC_1^M(A,\alpha_A)$ coincides with $HC^{\alpha}_1(A)$ when $(A,\alpha_A)$ is commutative.

\begin{theorem}\label{application}
Let $(A,\alpha_A)$ be a Hom-associative (non-commutative) algebra satisfying the $\alpha$-identity condition (\ref{condition}). Then there is an exact sequence of vector spaces
\begin{align*}
A \star \! HC_1^{\alpha}(A) &\to \Ker\big(A\!\star \!L^{\alpha}(A)\to L^{\alpha}(A)\big) \to \Ker\big(A\!\star \![A,A]\to [A,A]\big)\\
 &\to HC^{\alpha}_1(A)\to HC_1^M(A,\alpha_A)\to \frac{[A,A]}{[A,[A,A]]}\to 0.
\end{align*}
\end{theorem}
{\it Proof.}
By using Lemma \ref{lema 6.4} and Proposition \ref{exact-tensor-1} we have the commutative diagram of Hom-Lie algebras (written without $\alpha$ endomorphisms)

\[
 \xymatrix{
  & A\star HC_1^{\alpha}(A) \ar[r]  \ar[d]^{\psi_1}& A \star L^{\alpha}(A) \ar[r] \ar[d]^{\psi_2} & A\star [A,A] \ar[r]  \ar[d]^{{\psi_3}}& 0\\
0 \ar[r] & HC_1^{\alpha}(A) \ar[r] & L^{\alpha}(A) \ar[r]^{\psi} & [A,A] \ar[r] & 0.
}
\]
 Since $\Coker(\psi_3)= [A,A]/[A,[A,A]]$, $\Coker(\psi_2)=HC_1^M(A,\alpha_A)$, $\Coker(\psi_1)=HC_1^{\alpha}(A)$
and $\Ker(\psi_1)=A\star HC_1^{\alpha}(A)$, the assertion is a consequence of the Snake Lemma. \rdg
\medskip

Let us remark that if $\alpha_A$ is an epimorphism, then the term $A \star HC_1^{\alpha}(A)$ in the exact sequence of Theorem \ref{application} can be replaced by $A/[A,A] \otimes \! HC_1^{\alpha}(A)$ since they are isomorphic by Proposition \ref{proposition 5.10}. In particular, if $\alpha_A=\id_A$, the exact sequence in Theorem \ref{application} coincides with that of \cite[Theorem 5.7]{Gu}.
\newpage

\centerline{\bf Acknowledgements}

First and second authors were supported by Ministerio de Economía y Competitividad (Spain) (European FEDER support included), grant MTM2013-43687-P.   Second author was supported by Xunta de Galicia, grants EM2013/016 and GRC2013-045 (European FEDER support included) and by Shota Rustaveli National Science Foundation, grant DI/12/5-103/11.

\begin{center}

\end{center}

\end{document}